\documentclass[12pt,a4paper,reqno]{amsart}

\def\schur{\operatorname{{\sf Schur}}}

\usepackage{amsmath}
\usepackage{amssymb}
\usepackage{amsfonts}

\advance\textheight by \topskip
\title[Schur's old determinant]
{Schur's old determinant proves a brand--new
theorem of Garrett--Ismail--Stanton}

\author{Helmut Prodinger}

\address{ Helmut Prodinger,
The John Knopfmacher
Centre for Applicable Analysis and Number Theory,
 Department of Mathematics,
University of the Witwatersrand, P.~O. Wits, 
2050 Johannesburg, South Africa, email:
{\tt helmut@gauss.cam.wits.ac.za},\newline
homepage: {\tt http://www.wits.ac.za/helmut/index.htm}
}

\date{April 2, 2000}

\keywords{Rogers--Ramanujan identities, Schur's determinant}
\subjclass{05A30}

\begin{document}
\begin{abstract}
Garrett, Ismail, and Stanton gave a general formula that
contains the Rogers--Ramanjuan identities as special 
cases. We show how easy this is when using a determinant
that Schur introduced in 1917.
\end{abstract}

\maketitle

In a recent difficult paper
\cite{GaIsSt99} Garrett, Ismail, and Stanton
prove, amongst many other things, the following generalization
of the celebrated Rogers--Ramanujan identities:
\begin{align}\label{stanton}\begin{split}   
\sum_{n\ge0}\frac{q^{n^2+mn}}{(1-q)(1-q^2)\dots(1-q^n)}&
=(-1)^mq^{-\binom{m}{2}}E_{m-2}\prod_{n\ge0}
\frac{1}{(1-q^{5n+1})(1-q^{5n+4})}
\\&-(-1)^mq^{-\binom{m}{2}}D_{m-2}\prod_{n\ge0}
\frac{1}{(1-q^{5n+2})(1-q^{5n+3})},
\end{split}\end{align}
with the {\it Schur\/} polynomials, defined by
\begin{align*}
D_{m}&=D_{m-1}+q^mD_{m-2},  && D_0=1, \ D_1=1+q,\\
E_{m}&=E_{m-1}+q^mE_{m-2},  && E_0=1, \ E_1=1.\\
\end{align*}
(Another proof, based on generalized Engel expansions, 
can be found in \cite{AnKnPa99}.)

Schur \cite{Schur17} has computed the limits
\begin{equation*}   
D_\infty=\prod_{n\ge0}
\frac{1}{(1-q^{5n+1})(1-q^{5n+4})},
\qquad
E_\infty=\prod_{n\ge0}
\frac{1}{(1-q^{5n+2})(1-q^{5n+3})}.
\end{equation*}

The aim of this note is to give a simple proof of this
result, based on the following determinant of Schur:

\begin{equation*}   
\schur(x):=\left|
\begin{matrix}   
1&xq^{1+m}&&&&\dots\\
-1&1&xq^{2+m}&&&\dots\\
&-1&1&xq^{3+m}&&\dots\\
&&-1&1&xq^{4+m}&\dots\\
&&&\ddots&\ddots&\ddots
\end{matrix}
\right|.
\end{equation*}

Expanding the determinant with respect to the first column
(``top--recursion'') we get
\begin{equation*}   
\schur(x)=\schur(xq)+xq^{1+m}\schur(xq^2).
\end{equation*}
Setting
\begin{equation*}   
\schur(x)=\sum_{n\ge0}a_nx^n, 
\end{equation*}
we get, upon comparing coefficients, 
\begin{equation*}   
a_n=q^na_n+q^{1+m}q^{2n-2}a_{n-1},
\end{equation*}
or
\begin{equation*}   
a_n=\frac{q^{2n-1+m}}{1-q^n}a_{n-1}.
\end{equation*}
Since $a_0=1$, iteration leads to
\begin{equation*}   
a_n=\frac{q^{n^2+mn}}{(1-q)(1-q^2)\dots(1-q^n)},
\end{equation*} 
and thus
the left hand side of (\ref{stanton})
can be expressed by $\schur(1)$. 

On the other hand, $\schur(1)$ is the limit of the 
{\it finite\/}
determinants
\begin{equation*}   
\schur_n:=\left|
\begin{matrix}   
1&q^{1+m}&&&&\dots\\
-1&1&q^{1+m}&&&\dots\\
&-1&1&q^{2+m}&&\dots\\
&&-1&1&q^{3+m}&\dots\\
\vdots&\vdots&\vdots&\ddots&\ddots&\ddots\\
&&&-1&1&q^{n+m}\\
&&&&-1&1\\
\end{matrix}
\right|.
\end{equation*}
Expanding this determinant with respect to the last
row 
(``bottom--recursion'') we get
\begin{equation*}   
\schur_n=\schur_{n-1}+q^{n+m}\schur_{n-2}.
\end{equation*}
We see that the sequences ${\langle D_{n+m}\rangle}_n$
and ${\langle E_{n+m}\rangle}_n$ satisfy this recursion, 
and thus any linear combination; set
\begin{equation*}   
\schur_n=\lambda D_{n+m}+\mu E_{n+m}.
\end{equation*}
We can determine the parameters $\lambda$ and $\mu$
using the initial conditions $\schur_0=1$,  
$\schur_1=1+q^{1+m}$, which leads to
\begin{equation*}   
\lambda=\frac{E_m-E_{m-1}}{D_{m-1}E_m-D_mE_{m-1}}=
\frac{q^mE_{m-2}}{D_{m-1}E_m-D_mE_{m-1}}
\end{equation*}
and
\begin{equation*}   
\mu=\frac{D_m-D_{m-1}}{D_{m}E_{m-1}-D_{m-1}E_{m}}=
\frac{q^mD_{m-2}}{D_{m}E_{m-1}-D_{m-1}E_{m}}.
\end{equation*}
We can find the nicer form
\begin{equation*}   
\lambda=
\frac{q^mE_{m-2}}{D_{m-1}E_m-D_mE_{m-1}}=
(-1)^mq^{-\binom{m}{2}}E_{m-2},
\end{equation*}
which amounts to prove that
\begin{equation*}   
{D_{m-1}E_m-D_mE_{m-1}}=(-1)^mq^{\binom{m+1}{2}}.
\end{equation*}
This however is an easy induction on $m$, the beginning
$m=0$ being trivial. The induction step goes like this:
\begin{align*}
D_{m}E_{m+1}-D_{m+1}E_{m}&=
D_{m}(E_{m}+q^{m+1}E_{m-1})-(D_{m}+q^{m+1}D_{m-1})E_{m}\\
&=q^{m+1}(D_{m}E_{m-1}-D_{m-1}E_{m})\\
&=-q^{m+1}(-1)^mq^{\binom{m+1}{2}}=(-1)^{m+1}q^{\binom{m+2}{2}}.
\end{align*}  
Analogously we obtain
\begin{equation*}   
\mu=
\frac{q^mD_{m-2}}{D_{m}E_{m-1}-D_{m-1}E_{m}}=
(-1)^{m+1}q^{-\binom{m}{2}}D_{m-2}.
\end{equation*}
Hence
\begin{align*}
\schur_n=(-1)^mq^{-\binom{m}{2}}E_{m-2}D_{n+m}-
(-1)^mq^{-\binom{m}{2}}D_{m-2}E_{n+m}.
\end{align*}
Performing the limit $n\to\infty$ this turns into
\begin{align*}
\schur(1)=(-1)^mq^{-\binom{m}{2}}E_{m-2}D_{\infty}-
(-1)^mq^{-\binom{m}{2}}D_{m-2}E_{\infty},
\end{align*}
 which is (\ref{stanton}).

\bibliographystyle{plain}

\end{document}